\documentclass{crm-article}

\usepackage{citehack}
\usepackage{algorithm}
\usepackage{algorithmic}
\floatname{algorithm}{Алгоритм}

\def\ai#1{{\color{black}#1}} 
\def\vm#1{{\color{black}#1}} 

\newcommand{\argmin}{\mathop{\mathrm{argmin}}}
\newcommand{\argmax}{\mathop{\mathrm{argmax}}}

\newcommand{\RR}{\mathbb{R}}

\def\comments#1{{\color{black}#1}}

\begin{document}

\journalVol{10}

\journalNo{1}
\setcounter{page}{1}

\journalSection{Математические основы и численные методы моделирования}
\journalSectionEn{Mathematical modeling and numerical simulation}

\journalReceived{01.06.2016.}

\journalAccepted{01.06.2016.}

\UDC{519.85}
\title{\vm{Калибровка параметров модели расчета матрицы корреспонденций для г. Москвы.}}
\titleeng{\vm{Calibration of model parameters for calculating correspondence matrix for Moscow.}}


 \author[1,2]{\firstname{A.\,С.}~\surname{Иванова}}
 \authorfull{Анастасия Сергеевна Иванова}
 \authoreng{\firstname{A.\,S.}~\surname{Ivanova}}
 \authorfulleng{Anastasiya S. Ivanova}
 \email{anastasiya.s.ivanova@phystech.edu}

 \author[2]{\firstname{С.\,С.}~\surname{Омельченко}}
 \authorfull{Сергей Сергеевич Омельченко}
 \authoreng{\firstname{S.\,S.}~\surname{Omelchenko}}
 \authorfulleng{Sergey S. Omelchenko}
 \email{sergey.omelchenko@phystech.edu}
 \affiliation[1]{Национальный исследовательский университет «Высшая школа экономики»,\protect\\ Россия, 101000, г. Москва, адрес	Мясницкая улица, д. 20}
 \affiliationeng{Higher School of Economics — National Research University,\protect\\ 20 Myasnitskaya street, Moscow, 101000, Russia}

 \author[2]{\firstname{Е.\,В.}~\surname{Котлярова}}
 \authorfull{Екатерина Владимировна Котлярова}
 \authoreng{\firstname{E.\,V.}~\surname{Kotliarova}}
 \authorfulleng{Ekaterina V. Kotliarova}
 \email{kotlyarova.ev@phystech.edu}

 \author[2]{\firstname{В.\,В.}~\surname{Матюхин}}
 \authorfull{Владислав Вячеславович Матюхин}
 \authoreng{\firstname{V.\,V.}~\surname{Matyukhin}}
 \authorfulleng{Vladislav V. Matyukhin}
 \email{vladmatyukh@gmail.com}

 \affiliation[2]{Национальный исследовательский университет «Московский физико-технический институт»,\protect\\ Россия, 141701, г. Долгопрудный, Институтский пер., д. 9}
 \affiliationeng{National Research University Moscow Institute of Physics and Technology,\protect\\ 9 Institute lane, Dolgoprudny, 141701, Russia}

\thanks{Исследование выполнено при поддержкее Российского фонда фундаментальных исследований (проект 18-29-03071 мк). Исследование В.В. Матюхина выполнено при поддержке Министерства науки и высшего образования Российской Федерации (госзадание) №075-00337-20-03, номер проекта 0714-2020-0005.}

\thankseng{The research was supported by Russian Foundation for Basic Research (project 18-29-03071 mk). The research of V. V. Matyukhin was supported by the Ministry of Science and Higher Education of the Russian Federation (Goszadaniye) №075-00337-20-03, project No. 0714-2020-0005.}

\begin{abstract}
В данной работе рассматривается задача восстановления матрицы корреспонденци\vm{й} для наблюдений реальных корреспонденций в г. Москве. Следуя общепринятому подходу~\cite{gasnikov2013book}, транспортная сеть рассматривается как ориентированный граф, дуги которого соответствуют участкам дороги, а вершины графа - районы, из которых выезжают / в которые въезжают участники движения. Число жителей города считается постоянным.
Задача восстановления матрицы корреспонденци\vm{й} состоит в расчете всех корреспонденций из района $i$ в  район $j$.
 
Для восстановления матрицы предлагается использовать один из наиболее популярных в урбанистике способов расчета матрицы корреспонценци\vm{й} - энтропийная модель. В работе, базируясь на работе~\cite{Wilson1978}, приводится описание эволюционного обоснования энтропийной  модели, описывается основная идея перехода к решению задачи энтропийно-линейного программирования (ЭЛП) при расчете матрицы корреспонденци\vm{й}.  Для решения полученной задачи ЭЛП предлагается перейти к двойственной задач\vm{е} и решать задачу относительно двойственных переменных.  В работе описывается несколько численных методов оптимизации для решения данной задачи: \vm{алгоритм Синхорна} и ускоренный \vm{алгоритм Синхорна}. Далее приводятся численные эксперименты для следующих вариантов функций затрат: линейная функция затрат и \vm{сумма} степенной и логарифмической функции затрат. В данных функциях затраты представляют из себя некоторую комбинацию среднего времени в пути и расстояния между районами, которая зависит от параметров. Для каждого набора параметров функции затрат рассчитывается матрица корреспонденций и далее оценивается качество восстановленной матрицы относительно известной матрицы корреспонденций.
Мы предполагаем, что шум в восстановленной матрице корреспонденци\vm{й} является гауссовским, в результате в качестве метрики качества выступает среднеквадратичное отклонение. Данная задача представляет из себя задачу невыпуклой оптимизации. В статье приводится обзор безградиенных методов оптимизации для решения невыпуклых задач. Так как число параметров функции затрат небольшое, для  определение оптимальных параметров функции затрат было выбрано использовать метод перебора по сетке значений. Таким образом, для каждого набора параметров рассчитывается матрица корреспонденций и далее оценивается качество восстановленной матрицы относительно известной матрицы корреспонденций.  Далее по минимальному значению невязки для каждой функции затрат определяется для какой функции затрат и при каких  значениях параметров восстановленная матрица наилучшим образом описывает реальные корреспонденци\vm{й}.

\end{abstract}

\begin{abstracteng}
In this paper, we consider the problem of restoring the correspondence matrix based on the observations of real correspondences in Moscow. Following the conventional approach~\cite{gasnikov2013book}, the transport network is considered as a directed graph whose edges correspond to road sections and the graph vertices correspond to areas that the traffic participants leave or enter. The number of city residents is considered constant. 
The problem of restoring the correspondence matrix is to calculate all the correspondence from the $i$ area to the $j$ area. 

To restore the matrix, we propose to use one of the most popular methods of calculating the correspondence matrix in urban studies - the entropy model. In our work, which is based on the work~\cite{Wilson1978}, we describe the evolutionary justification of the entropy model and the main idea of the transition to solving the problem of entropy-linear programming (ELP) in calculating the correspondence matrix. To solve the ELP problem, it is proposed to pass to the dual problem. In this paper, we describe several numerical optimization methods for solving this problem: the Sinkhorn method and the Accelerated Sinkhorn method. We provide numerical experiments for the following variants of cost functions: a linear cost function and a superposition of the power and logarithmic cost functions. In these functions, the cost is a combination of average time and distance between areas, which depends on the parameters. The correspondence matrix is calculated for multiple sets of parameters and then we calculate the quality of the restored matrix relative to the known correspondence matrix. 

We assume that the noise in the restored correspondence matrix is Gaussian, as a result, we use the standard deviation as a quality metric. 
The article provides an overview of gradient-free optimization methods for solving non-convex problems. Since the number of parameters of the cost function is small, we use the grid search method to find the optimal parameters of the cost function. Thus, the correspondence matrix calculated for each set of parameters and then the quality of the restored matrix is evaluated relative to the known correspondence matrix. Further, according to the minimum residual value for each cost function, we determine for which cost function and at what parameter values the restored matrix best describes real correspondence.
\end{abstracteng} 

\keywordeng{
Correspondence matrix calculation model}
\keywordeng{Entropy Linear Programming}
\keywordeng{Sinkhorn method}
\keywordeng{Accelerated Sinkhorn method}

\maketitle

\paragraph{Введение}

 Методы моделирования матриц корреспонденций начали активно развиваться в 60-х годах прошлого века. Задача состояла в том, чтобы по $2n$ параметрам (численности активного населения в $n$ районах и рабочих мест в этих районах) определить $n^{2}$ параметров - матрицу корреспонденций. Модели расчета матрицы корреспонденций базируются на характеристиках районов ($2n$ параметров) и \vm{матрицах затрат} ($n^{2}$ известных чисел, характеризующих \vm{затраты на дорогу из одного района в другой: время, расстояние и т.д.}). Самый первой моделью являлась гравитационная модель, в основу которой был положен аналог закона всемирного тяготения  Ньютона-Гука. Однако, в отличие от отмеченного закона физики, параметры этого закона (показатель степени, с которым входят затраты между районами в оценку корреспонденци\vm{й}) неизвестны. В данной работе рассматривается один из наиболее популярных в урбанистике способов расчета матрицы корреспонценци\vm{й}, имеющий более чем сорокалетнюю историю, - энтропийная модель. В основе рассматриваемой модели лежит марковский процесс в пространстве огромной размерности, то есть данный процесс порождает макросистему. Данный марковский процесс представляет собой ветвящийся процесс специального вида: "модель стохастической химической кинетики". Вывод энтропийной модели приведен на базе бинарных реакций обменного типа, популярных в различного рода физических и социально-экономических приложениях моделей стохастической химической кинетики~\cite{Gardiner2009, Weidlich2010}. 
 В данной статье на базе энтропийных моделей расчета матрицы корреспонден\vm{ций} (см. книгу А. Дж. Вильсона ~\cite{Wilson1978}, книгу Е.В. Гасниковой и А.В. Гасникова ~\cite{gasnikov2013book}) предложен способ определения параметров таких моделей и выбор наилучшей. Стоит отметить, что ранее такого рода исследования сосредотачивались в основном только на гравитационной и энтропийной модели, у которых имеется всего один параметр для калибровки. \vm{Данная работа сосредоточена на более общих и точных современных вариантах данных моделей}, описанных в цикле работ Е.В. Гасниковой и А.В. Гасникова (~\cite{rangePageRank}, ~\cite{gasnikov2013book}). \vm{Калибровка таких моделей на реальных данных (были использованы данные по районам Москвы) ранее не осуществлялась}. Отметим, что при описании процедуры калибровки для получения гарантированного (теоретически обоснованного) результата необходим аккуратный анализ возникающих задач оптимизации и используемых алгоритмов. Предложенный в статье анализ базируется на современных достижениях в двух направлениях: "Оптимальный транспорт" (современный анализ оценок скорости сходимости алгоритма Синхорна, ~\cite{ComputationalOptimalTransport2018}) и безградиентные методы невыпуклой оптимизации с неточным оракулом (\cite{GradientFreeOptimizationMethods}). Отметим также метод балансировки~\cite{Shele, Bregman}, который эквивалентен алгоритму Синхорна и может рассматриваться для решения задачи ЭЛП. 

Статья организована следующим образом.

В п.2 описан вывод энтропийной модели. В результате, показано, что задачу расчета матрицы корреспонденци\vm{й} можно рассматривать как задачу энтропийно-линейного программирования. 

В п.3 описано решение задачи энтропийно-линейного программирования из п.2. Предлагается перейти к двойственной задач\vm{е}, которую предлагается решать \vm{Алгоритмом Синхорна}, попеременно минимизируя \vm{двойственную целевую функцию} по одному из двух блоков двойственных переменных. Также описан \vm{ускоренный алгоритм} Синхорна, который позволяет быстрее решать рассматриваемую задачу.

В п.4 описывается задача оценки невязки между реальной матрицей и матрицей, восстановленной в результате решения задачи из п.3. 

В п.5 проведены численные эксперименты по расчету матрицы корреспонденци\vm{й} по реальным данным для города Москвы для разных функций затрат. Проводится сравнительный анализ для определения функции затрат, которая наилучшим образом описывает реальные данные.  

\paragraph{Постановка задачи. Энтропийная модель расчета матрицы корреспонденци\vm{й}}
В данном разделе, базируясь на работе~\cite{Wilson1978}, приведем эволюционное обоснование энтропийной модели. \vm{Дальнейшее изложение вывода идет согласно статье~\cite{gasnikov2016}.}

Пусть в некотором городе имеется $n$ районов. Общее число жителей города постоянно и равно $N$, при это выполняется $N \gg n^2$. Пусть $L_i \geq 0$ это число жителей, выезжающих в типичный день за рассматриваемый промежуток времени из района $i$, а $W_j \geq 0$ число жителей, приезжающих на работу в район $j$ в типичный день за рассматриваемый промежуток времени. В рамках рассматриваемого подхода данные величины являются входными параметрами для модели, т.е. они не моделируются.  При этом будут выполняться следующие соотношения: $\sum\limits_{i=1}^{n} L_i = \sum\limits_{j=1}^{n} W_j = N$.

Обозначим через $d_{ij}(t) \geq 0$ -- число жителей, живущих в $i$-м районе и работающих в $j$-м в момент времени $t$. Мы предполагаем, что со временем жители могут меняться только квартирами, поэтому во все моменты времени $t\geq 0$ выполнено 
\begin{eqnarray*}
   d_{ij}(t) \geq 0, \, \, \sum\limits_{i,j=1}^{n,n} d_{ij}(t)  \vm{=} N, \, \,  \sum\limits_{j=1}^{n} d_{ij}(t) \vm{=} L_i, \, \, \sum\limits_{i=1}^{n} d_{ij}(t) \vm{=} W_j, \, \, i, j = 1, \ldots, n.
\end{eqnarray*}
Определим следующее множество:
\begin{eqnarray}
\label{Q_set}
   Q = \left\{ d_{ij} \geq 0: \, \, \sum\limits_{i,j=1}^{n,n} d_{ij} \vm{=} N, \, \, \sum\limits_{j=1}^{n} d_{ij} \vm{=} L_i, \, \, \sum\limits_{i=1}^{n} d_{ij} \vm{=} W_j, \, \, i, j = 1, \ldots, n \right\}.
\end{eqnarray}

Отметим, что основным стимулом к обмену места жительства для жителя города будут являться транспортные издержки, то есть для каждого жителя работать далеко от дома плохо из-за больших транспортных издержек. Будем считать, что эффективной функцией затрат~\cite{gasnikov2013book}  является функция $R(T)= \tfrac{\gamma T}{2}$, где $T > 0$ - затраты на путь от дома до работы, которые определяются как временем так и расстоянием, а $\gamma>0$ -- настраиваемый параметр модели, который можно интерпретировать как цену единицы затрат на путь от работы до дома. \vm{Далее в работе будем подразумевать под функцией затрат только $T(\alpha, \beta, \gamma)$, где  $\alpha, \beta, \gamma$ -- настраиваемые параметры модели.}

Динамику процесса можно описать следующим образом: пусть в момент времени $t \geq 0$ $r$-й житель живет в $k$-м районе и работает в $m$-м, а $s$-й житель живет в $p$-м районе и работает в $q$-м. Тогда $\lambda_{k,m;p,q}(t)\Delta t + o(\Delta t)$- есть вероятность того, что жители с номерами  $r$ и $s$ ($1\leq r < s \leq N$)  "поменяются" квартирами в промежутке времени $(t, t+ \Delta t)$.  Вероятность обмена местами жительства зависит только от мест проживания и работы обменивающихся:
\begin{eqnarray*}
    \lambda_{k,m;p,q}(t) = \lambda_{k,m;p,q} = \lambda N^{-1} \exp\Bigl(\underbrace{R(T_{km}) + R(T_{pq})}_{\substack{\text{cуммарные затраты} \\ \text{  до обмена}}} - \underbrace{(R(T_{pm}) + R(T_{kq}))}_{\substack{\text{cуммарные затраты} \\ \text{ после обмена}}} \Bigr) > 0,
\end{eqnarray*}
где коэффициент $0 < \lambda = O(1)$ характеризует интенсивность обменов. Отметим, что совершенно аналогичным образом можно было рассматривать случай, когда жители могут обмениваться местами работы. \comments{То есть мы предполагаем некое равноправие агентов (жителей) внутри фиксированной корреспонденций и их независимость  ~\cite{gasnikov2014, ethier1986}.} 

Согласно эргодической теореме для марковских цепей (в независимости от 
начальной конфигурации $\{ d_{ij}(0)\}_{i,j=1,1}^{n,n}$) ~\cite{Malyshev2008, gasnikov2013book, Gardiner2009, Weidlich2010, Sandholm2010, Borovkov1999, levin2009}
предельное распределение совпадает со стационарным (инвариантным), которое можно посчитать (получается
проекция прямого произведение распределений Пуассона на $Q$):
\begin{eqnarray*}
    \lim\limits_{t \rightarrow \infty} P(d_{ij}(t) = d_{ij}, \, i,j =1, \ldots, n)& =&Z^{-1} \prod\limits_{i,j=1,1}^{n,n}\exp{(-2R(T_{ij})d_{ij})}\cdot (d_{ij}!)^{-1} \\ &\overset{def}{=} &p(\{ d_{ij}\}_{i,j=1,1}^{n,n}),
\end{eqnarray*}
где $\{ d_{ij}\}_{i,j=1,1}^{n,n}\in Q$, а "статсумма"~$Z$. \comments{Отметим, что начальная конфигурация $\{ d_{ij}(0)\}_{i,j=1,1}^{n,n}$ влияет на время выхода на стационарное состояние.} При этом стационарное распределение
$p(\{ d_{ij}\}_{i,j=1,1}^{n,n})$ удовлетворяет условию детального равновесия~\cite{gasnikov2014, Sandholm2010}:
\begin{eqnarray*}
     (d_{km}+1)(d_{pq}+1)p\left(\left\{d_{11}, ..., d_{km}+1,...,d_{pq}+1,..., d_{pm}-1, ...,  d_{kq}-1, ..., d_{nn}  \right\} \right) \lambda_{k,m;p,q}\\ =d_{pm}d_{kq}p\left(\{ d_{ij}\}_{i,j=1,1}^{n,n} \right) \lambda_{k,m;p,q}\vm{.}
\end{eqnarray*}
При $N \gg 1$ распределение $p(\{ d_{ij}\}_{i,j=1,1}^{n,n})$ экспоненциально сконцентрировано на множестве $Q$ в  $O(\sqrt{N})$ окрестности наиболее вероятного значения $d^* = \{ d^*_{ij}\}_{i,j=1,1}^{n,n}$, которое определяется, как решение задачи энтропийно-линейного программирования (ЭЛП)~\cite{Malyshev2008, gasnikov2013book }:
\begin{eqnarray*}
    \ln p(\{ d_{ij}\}_{i,j=1,1}^{n,n}) \vm{=} \max\limits_{\{ d_{ij}\}_{i,j=1,1}^{n,n} \in Q} - \gamma \sum\limits_{i,j = 1,1}^{n,n}d_{ij}T_{ij} - \sum\limits_{i,j = 1,1}^{n,n}d_{ij}\ln d_{ij}.
\end{eqnarray*}
Это следует из теоремы Санова о больших уклонениях для мультиномиального распределения~\cite{Sanov1957}. 
Отметим, что в данном параграфе описывается основная идея перехода к решению задачи ЭЛП при расчете матрицы корреспонденци\vm{й}. Более формально о полученном результате можно найти в \cite{gasnikov2016}.

\paragraph{Методы для решения задачи ЭЛП}

В данном разделе приведем описание способов решения задачи ЭЛП, возникающей при расчете матрицы корреспонденци\vm{й}.

Как было показано в предыдущем разделе, задачу восстановления матрицы корреспонденци\vm{й} можно записать как следующую задачу оптимизации:
\begin{equation}
\label{prob_st_init}
    \min_{d_{ij} \in Q} f(d_{ij}) := \gamma \sum\limits_{i,j = 1,1}^{n,n}d_{ij}T_{ij} + \sum\limits_{i,j = 1,1}^{n,n}d_{ij}\ln d_{ij},
\end{equation}
\comments{где $Q$ определяется как \eqref{Q_set} и  $T_{ij} := T_{ij}(\eta)$ -- функция затрат на перемещение из района $i$ в район $j$, которая зависит от вектора параметров $\eta$}. 

Введем следующую нормировку: $\sum\limits_{i,j = 1,1}^{n,n}d_{ij} = 1$, тогда ограничения можно переписать в следующем виде  $\sum\limits_{j=1}^n d_{ij} = l_i$ и $\sum\limits_{i=1}^n d_{ij} = w_j$, где $l_i = \tfrac{L_i}{N}$ $w_j = \tfrac{W_j}{N}$. И определим следующее множество 
\begin{eqnarray*}
   \tilde{ Q} = \left\{ d_{ij} \geq 0: \, \, \sum\limits_{i,j=1}^{n,n} d_{ij} \vm{=} 1, \, \, \sum\limits_{j=1}^{n} d_{ij} \vm{=} l_i, \, \, \sum\limits_{i=1}^{n} d_{ij} \vm{=} w_j, \, \, i, j = 1, \ldots, n \right\}.
\end{eqnarray*}
\comments{После введения нормировки, получаем что для задачи~\eqref{prob_st_init}:
\begin{eqnarray*}
    \gamma \sum\limits_{i,j = 1,1}^{n,n}N \cdot d_{ij}T_{ij} &+ & \sum\limits_{i,j = 1,1}^{n,n}N \cdot d_{ij}\ln N \cdot d_{ij} \\ &=& N \Bigl\{ \gamma \sum\limits_{i,j = 1,1}^{n,n}  d_{ij}T_{ij} +  \sum\limits_{i,j = 1,1}^{n,n} d_{ij}\ln d_{ij} +  \underbrace{\sum\limits_{i,j = 1,1}^{n,n} d_{ij}}_{=1}\ln N \Bigr\}.
\end{eqnarray*}
}
Тогда задача~\eqref{prob_st_init} перепишется в следующем эквивалентном виде:
\begin{equation}
\label{prob_st}
    \min_{ d_{ij} \in \tilde{ Q}} \gamma \sum\limits_{i,j = 1,1}^{n,n}d_{ij}T_{ij} +  \sum\limits_{i,j = 1,1}^{n,n}d_{ij}\ln d_{ij}.
\end{equation}
Далее введем два блока двойственных переменных $\lambda^{l} \in \RR^n$ и $\lambda^{w} \in \RR^n$, где $\lambda_i^{l}$ множитель к
 ограничению $\sum\limits_{j=1}^n d_{ij} = l_i$ и $\lambda_j^{w}$ множитель к ограничению $\sum\limits_{i=1}^n d_{ij} = w_j$. Применим для решения задачи~\eqref{prob_st} метод множителей Лагранжа. Для этого запишем двойственную задачу:

\begin{align*}
 \min_{d_{ij} \in \tilde{ Q}} \gamma \sum\limits_{i,j = 1,1}^{n,n}d_{ij}T_{ij} +  \sum\limits_{i,j = 1,1}^{n,n}d_{ij}\ln d_{ij}  \\ = \max\limits_{\lambda^{l}, \, \lambda^{w}} \min_{\substack{ \sum\limits_{i,j = 1,1}^{n,n}d_{ij} = 1,\\  d_{ij} \geq 0 }}   \Bigl\{ \gamma \sum\limits_{i,j = 1,1}^{n,n}d_{ij}T_{ij} +  \sum\limits_{i,j = 1,1}^{n,n} d_{ij}\ln d_{ij} + \sum\limits_{i = 1}^{n} \lambda^l_i (\sum\limits_{j = 1}^{n} d_{ij} - l_i)+ \sum\limits_{j = 1}^{n}  \lambda^w_j (\sum\limits_{i = 1}^{n} d_{ij} - w_j)   \Bigr\} \\ =
  \max\limits_{\lambda^{l}, \, \lambda^{w}} \Bigl\{- \langle \lambda^l, l\rangle -  \langle
\lambda^w, w\rangle  +  \min_{d_{ij} \geq 0} \Bigl\{  \sum\limits_{i,j = 1,1}^{n,n}d_{ij} \left( \gamma T_{ij}  +   \ln d_{ij} + \lambda^l_i + \lambda^w_j\right) + \nu(\sum\limits_{i,j = 1,1}^{n,n} d_{ij} - 1)   \Bigr\}  \Bigr\} \\ = 
\max\limits_{\lambda^{l}, \, \lambda^{w}} \Bigl\{- \langle \lambda^l, l\rangle -  \langle
\lambda^w, w\rangle  +   \Bigl\{  \sum\limits_{i,j = 1,1}^{n,n}d_{ij}( \lambda^l , \lambda^w ) \left( \gamma T_{ij}  +   \ln d_{ij}( \lambda^l , \lambda^w ) + \lambda^l_i + \lambda^w_j + \nu \right) - \nu   \Bigr\}  \Bigr\},
\end{align*}
где 
\begin{eqnarray*}
    d_{ij}( \lambda^l , \lambda^w ) = \argmin_{d_{ij} \geq 0} \Bigl\{  \sum\limits_{i,j = 1,1}^{n,n}d_{ij} \left( \gamma T_{ij}  +   \ln d_{ij} + \lambda^l_i + \lambda^w_j\right) + \nu(\sum\limits_{i,j = 1,1}^{n,n} d_{ij} - 1)   \Bigr\}  \Bigr\}.
\end{eqnarray*}
Используя условия оптимальности, получаем
\begin{eqnarray*}
    \gamma T_{ij} + \ln d_{ij} + \lambda_i^l + \lambda_j^w + 1 + \nu = 0, \quad \sum\limits_{i,j = 1,1}^{n,n} d_{ij} = 1.
\end{eqnarray*}
Решая данную систему уравнений и переопределяя $\lambda^l := - \lambda^l - \tfrac{1}{2}$ и $\lambda^w := - \lambda^w - \tfrac{1}{2}$ получаем, что 
\begin{eqnarray*}
    d_{ij}( \lambda^l , \lambda^w ) = \tfrac{ \exp(- \gamma T_{ij} + \lambda_i^l + \lambda_j^w )}{\sum\limits_{i,j = 1,1}^{n,n} \exp(- \gamma T_{ij} + \lambda_i^l + \lambda_j^w )} = \tfrac{ B_{ij}(\lambda^l , \lambda^w) }{\mathbf{1}^{T}B(\lambda^l , \lambda^w )\mathbf{1} } ,
\end{eqnarray*} 
где $B_{ij}(\lambda^l , \lambda^w ) = \exp(- \gamma T_{ij} + \lambda_i^l + \lambda_j^w )$. 
Подставляя это в двойственную задачу, получаем, что двойственная задача имеет вид 
\begin{eqnarray*}
    \max\limits_{\lambda^{l}, \, \lambda^{w}} \tilde \varphi ( \lambda^l , \lambda^w ) :=  \langle \lambda^l, l\rangle +  \langle
\lambda^w, w\rangle  - \ln \left(\mathbf{1}^{T}B(\lambda^l , \lambda^w )\mathbf{1}  \right)
\end{eqnarray*}
Перепишем задачу, как задачу минимизации с точностью до знака
\begin{eqnarray}
\label{dual_prob}
    \min\limits_{\lambda^{l}, \, \lambda^{w}}  \varphi ( \lambda^l , \lambda^w ) :=  \ln \left(\mathbf{1}^{T}B(\lambda^l , \lambda^w )\mathbf{1}  \right)  - \langle \lambda^l, l\rangle -  \langle
\lambda^w, w\rangle.
\end{eqnarray}
Для решения двойственной задачи рассмотрим \vm{Метод Альтернированной Минимизации} (Алгоритм~\ref{alg:AM1}). Для удобства описания алгоритма введем следующее обозначение. Множество $\{1, . . . , n\}$ векторов ${e_i}_{i=1}^n$ ортонормированного базиса разделено на $p$ непересекающихся блоков $I_k , \, \, k \in \{1, \ldots , p\}$. Пусть $S_k (x) = x + \text{span}\{e_i : i \in I_k\}$, подпространство, содержащее $x$ построенное на базисных векторах $k$-го блока.

\begin{algorithm}
\caption{Метод \vm{Альтернированной Минимизации}}
\begin{algorithmic}[1]\label{alg:AM1}
\STATE \textbf{Input:} $x^0$ -- starting point.
\FOR{$k \geq 0$} 
\STATE Choose $i_k \in 1, \ldots, p$.
\STATE Compute $x^{k+1} = \argmin\limits_{x \in S_{i_k}(x^k)} f(x)$.
\ENDFOR
\STATE \textbf{Output:} $x^k$.
\end{algorithmic}
\end{algorithm}
Отметим, что основной идеей данного алгоритма является минимизация по произвольно выбранному блоку переменных $i_k$ на каждой итерации. Для задачи~\eqref{dual_prob} мы будет рассматривать минимизацию по двум блокам: $\lambda^l$ и $\lambda^w$. Согласно Лемме~5 из~\cite{guminov2020} шаг минимизации по блоку $\lambda^l$ можно представить в следующем виде
\begin{eqnarray*}
    [\lambda^l ]^{k+1} = [\lambda^l ]^{k} + \ln(l) - \ln\left(B([\lambda^l ]^{k} , [\lambda^w ]^{k} )\mathbf{1}   \right),
\end{eqnarray*}
аналогичным образом можно представить шаг минимизации по блоку $\lambda^w$:
\begin{eqnarray*}
    [\lambda^w ]^{k+1} = [\lambda^w ]^{k} + \ln(w) - \ln\left(B^{T}([\lambda^l ]^{k} , [\lambda^w ]^{k} )\mathbf{1}   \right).
\end{eqnarray*}
\vm{Учитывая это, для решения~\eqref{dual_prob} получаем Алгоритм~\ref{alg:AM}.}

\begin{algorithm}
\caption{\vm{Алгоритм Синхорна}}
\begin{algorithmic}[1]\label{alg:AM}
\STATE \textbf{Input:} $x^0 = [[\lambda^l ]^{0}, [\lambda^w ]^{0}] = (0, \ldots, 0) \in \RR^{2n}$ -- starting point.
\FOR{$k \geq 0$} 
\IF{$k \, \, \text{mod} \, \, 2 = 0 $}
\STATE Compute 
\begin{eqnarray*}
    & [\lambda^l ]^{k+1} & = [\lambda^l ]^{k} + \ln(l) - \ln\left(B([\lambda^l ]^{k} , [\lambda^w ]^{k} )\mathbf{1}   \right),  \\
    & [\lambda^w ]^{k+1} &=  [\lambda^w ]^{k}.
    \end{eqnarray*}
 \ELSE 
 \STATE Compute 
\begin{eqnarray*}
    & [\lambda^l ]^{k+1} &=   [\lambda^l ]^{k},  \\
    & [\lambda^w ]^{k+1} & = [\lambda^w ]^{k} + \ln(w) - \ln\left(B^{T}([\lambda^l ]^{k} , [\lambda^w ]^{k} )\mathbf{1}   \right).
    \end{eqnarray*}
\ENDIF

\ENDFOR
\STATE \textbf{Output:} $x^k = [[\lambda^l ]^{k}, [\lambda^w ]^{k}] \in \RR^{2n}$.
\end{algorithmic}
\end{algorithm}

Отметим, что Алгоритм \vm{Альтернированной Минимизации} для задачи~\eqref{dual_prob} является хорошо \vm{известный алгоритм Синхорна} ~\cite{marco2013}.

Также для оптимального решения задачи ЭЛП будем рассматривать ускоренный вариант метода \vm{Альтернативной Минимизации}. Согласно~~\cite{guminov2020}, в качестве основы ускоренного метода Альтернативной Минимизации используется традиционный адаптивный \vm{ускоренный градиентный метод}. Для этой задачи этот вариант метода оказался быстрее на практике, чем другие способы ускорения. Здесь мы не используем одномерную минимизацию, чтобы найти размер шага, а вместо этого мы адаптируемся к константе Липшица $L$. Анализ скорости сходимости этого алгоритма можно найти в ~\cite{guminov2020}. В нашем случае  Ускоренный Метод \vm{Альтернативной Минимизации} \vm{представлен в виде Алгоритма~\ref{alg:AAM}}.

\begin{algorithm}
\caption{Ускоренный \vm{алгоритм Синхорна}}
\begin{algorithmic}[1]\label{alg:AAM}
\STATE \textbf{Input:} $x^0 := [[x^l ]^{0}, [x^w ]^{0}] = (0, \ldots, 0) \in \RR^{2n}$ -- starting point, $L_{0} = 1$, $a_{0} = 0$.

\REPEAT 
\STATE Set $y^0 := [[y^l ]^{0}, [y^w ]^{0}] = x^0$.
\STATE Set $v^0 := [[v^l ]^{0}, [v^w ]^{0}] = x^0$.
\STATE $L_{k+1}=L_{k}/2$

\WHILE{True}
   \STATE Set $a_{k+1}=\frac{1}{2L_{k+1}}+\sqrt{\frac{1}{4L^2_{k+1}}+a_k^2\frac{L_k}{L_{k+1}}}$
   \STATE Set $\tau_k=\frac{1}{a_{k+1}L_{k+1}}$
   \STATE Set $y^{k}=\tau_k v^k+(1-\tau_k)x^k$
   \STATE Choose $i_k = \argmax\limits_{i\in\{1,2\}} \|\nabla_i \varphi(y^k)\|^2$
    \IF{$i_k = 1 $}
    \STATE Compute 
    \begin{eqnarray*}
        & [x^l ]^{k+1} & = [y^l ]^{k} + \ln(l) - \ln\left(B([y^l ]^{k} , [y^w ]^{k} )\mathbf{1}   \right),  \\
        & [x^w ]^{k+1} &=  [y^w ]^{k}.
        \end{eqnarray*}
     \ELSE 
     \STATE Compute 
    \begin{eqnarray*}
        & [x^l ]^{k+1} &=   [y^l ]^{k},  \\
        & [x^w ]^{k+1} & = [y^w ]^{k} + \ln(w) - \ln\left(B^{T}([y^l ]^{k} , [y^w ]^{k} )\mathbf{1}   \right).
        \end{eqnarray*}
    \ENDIF
    \STATE Set $v^{k+1} = v^{k} - a_{k+1}\nabla \varphi(y^k)$
    \IF{$\varphi(x^{k+1})\leqslant \varphi(y^{k})-\frac{\|\nabla \varphi (y^k)\|^2}{2L_{k+1}}$}
        \STATE Set $\hat{d}^{k+1} = \frac{a_{k+1}d^k(y^k)+L_{k}a_{k}^2\hat{d}^k}{L_{k+1}a_{k+1}^2}$
        \STATE \textbf{break}
    \ENDIF
    \STATE Set $L_{k+1}=2L_{k+1}$.
\ENDWHILE

\UNTIL{ $|f(\hat{d}^{k+1}) + \varphi(x^{k+1})| \leq \varepsilon_{f}$,  $||\hat{d}^{k+1} \mathbf{1} - l||_2 \leq \varepsilon_{eq}$, $||(\hat{d}^{k+1})^{T} \mathbf{1} - w||_2 \leq \varepsilon_{eq}$}

\STATE \textbf{Выход:} $\hat{d}^{k+1}$, $x^{k+1}$.
\end{algorithmic}
\end{algorithm}

При этом вектор градиента \comments{функции \eqref{dual_prob}} представляет из себя следующий вектор: 
\begin{eqnarray*}
    && \nabla \varphi (\lambda^l , \lambda^w) = \left[\nabla_{1} \varphi^{T}, \nabla_{2} \varphi^{T} \right]^{T}, \\ && \text{где}  \, \, \,  \nabla_{1} \varphi(\lambda^l , \lambda^w) = - l + \frac{B(\lambda^l , \lambda^w)\mathbf{1}}{\mathbf{1}^{T}B(\lambda^l , \lambda^w)\mathbf{1}},\quad 
    \nabla_{2} \varphi (\lambda^l , \lambda^w)= - w + \frac{B^{T}(\lambda^l , \lambda^w)\mathbf{1}}{\mathbf{1}^{T}B(\lambda^l , \lambda^w)\mathbf{1}}.
\end{eqnarray*}

\paragraph{Задача подсчёта невязки для восстановленной матрицы затрат}

В данном параграфе опишем задачу подсчета невязки между восстановленной матрицей по затратам и реальной матрицей корреспонденци\vm{й}. Подсчет невязки необходим, чтобы оценить насколько хорошо выбранная функция затрат описывает реальные данные.

Для постановки задачи подсчёта невязки между $d_{ij}$ -- исходной матрицей корреспонденций и  $\widehat{d}_{ij}(\alpha)$ -- восстановленной матрицей корреспонденций, \comments{предположим, что в восстановленной матрице корреспонденций  $\widehat{d}_{ij}(\alpha)$ шум является гауссовским (мы восстанавливаем матрицу неточно, с шумом). Тогда восстановленную матрицу можно рассматривать как нормально распределенную выборку $N(\theta, \sigma^2)$,} и плотность вероятности нормального распределения можно рассчитать следующим образом:
\comments{
\begin{eqnarray*}
    &p(x) = \dfrac{1}{ \sqrt{2\pi \sigma^2}} \cdot e^{\tfrac{-(x - \theta)^2}{2 \sigma^2}  },
\end{eqnarray*}
где в качестве матожидания введенного гауссовского распределения будет выступать исходная матрица корреспонденций $d_{ij}$.}

Следовательно, максимизируя правдоподобие, получим:
\begin{eqnarray*}
    L(\widehat{d}_{ij}(\alpha)) & = & \prod_{i,j=1}^{n,n} (p(d_{ij}(\alpha))) = \prod_{i,j = 1}^{n,n} {\frac{1}{\sqrt{2\pi\sigma^2}} \cdot e^{-(d_{ij} - \widehat{d}_{ij}(\alpha))^2/2\comments{\sigma^2}  }} \rightarrow \max \\
    \log{L(\widehat{d}_{ij}(\alpha))} &=& \sum_{i,j=1}^{n,n} \left\{-\log{\sqrt{2\pi\sigma^2}} - \frac{1}{2\comments{\sigma^2}} \cdot (d_{ij} - \widehat{d}_{ij}(\alpha))^2 \right\}\rightarrow \max    \\
    \log{L(\widehat{d}_{ij}(\alpha))} &=& -\log{\sqrt{2\pi\sigma^2}} - \frac{1}{2\comments{\sigma^2}}  \cdot \sum_{k=1}^{n}{ (d_{ij} - \widehat{d}_{ij}(\alpha))^2 \rightarrow \max  } \\
\end{eqnarray*}
Изменив знак перед выражением, получаем следующую задачу для подсчёта невязки:
\begin{eqnarray}
    \min\limits_{\alpha}\sum\limits_{i, j=1, 1}^{n, n}(d_{ij} - \widehat{d}_{ij}(\alpha))^2 \, \,    \text{или в нормированном случае} \, \, \frac{ \min\limits_{\alpha}\sum\limits_{i,j=1,1}^{n,n}(d_{ij} - \widehat{d}_{ij}(\alpha))^2}{n^2}.
    \label{error_of_restored_matrix}
\end{eqnarray}
Отметим что данная задача является задачей минимизации, которая зависит от параметра $\alpha$ (это может быть вектор параметров, в зависимости от количества параметров в рассматриваемой функции затрат). 

Целевая функция полученной задачи является невыпуклой функцией. Для решения данной задачи предлагается использовать безградиентные методы.
В частности в рассматриваемой задач\vm{е} для поиска оптимального параметра (параметров) $\alpha$ в возникающей при подсчёте невязки задаче минимизации, используется метод перебора так как число параметров в зависимости от функции затрат $1-3$. 

Однако, в качестве обзора, приведем описание еще нескольких безградиентных методов для задач невыпуклой оптимизации.
Рассмотрим метод имитации отжига, для работы которого не требуется гладкость функции ~\cite{AnnealingMultistart}. Он является вариантом метода случайного поиска и известен как алгоритм Метрополиса. Для задач оптимизации имитация процесса может быть произведена следующим образом. Вводится параметр $T$, который имеет смысл температуры, и в начальный момент ему устанавливается значение $T_0$. Набор переменных, по которым происходит оптимизация, будет обозначаться как $x$. В качестве начального состояния системы выбирается произвольная точка. Далее запускается итерационный процесс — на каждом шаге из множества соседних состояний случайно выбирается новое $\vm{\hat{x}}$. Если значение функции в этой точке меньше, чем значение в текущей точке, то эта точка выбирается в качестве нового состояния системы. В ином случае (т.е. если $\vm{f(\hat{x})} > f(x)$) такой переход происходит с вероятностью $P$, зависящей от температуры $T$, текущего состояния и кандидата на новое состояние $\vm{\hat{x}}$ следующим образом:

\begin{eqnarray*}
    P = e^{-\frac{\vm{f(\hat{x})} - f(x)}{T}}.
\end{eqnarray*}

Также стоит упомянуть метод ломаных, который применим к классу функций одной переменной, удовлетворяющих условию Липшица ~\cite{BrokenCurve}. Говорят, что функция $f(x)$ удовлетворяет условию Липшица, если найдётся такая константа $L > 0$, что:
\begin{eqnarray*}
    \mid f(x) - f(y)\mid \leq L \cdot \mid x - y\mid \quad \forall x, y \in [a, b]
\end{eqnarray*}
Пусть функция $f(x)$ удовлетворяет условию Липшица на отрезке $[a, b]$. Зафиксируем какую-либо точку $y \in [a, b]$ и определим функцию $g(x, y) = f(y) - L \cdot \mid x - y\mid$ переменной $a \leq x \leq b$. Функция $g(x, y)$ кусочно-линейна на $[a, b]$, и график её представляет ломаную линию, составленную из отрезков двух прямых, имеющих угловые коэффициенты $L$ и $-L$ и пересекающихся в точке $(y, f(y))$. Также в силу липшицевого условия:
\begin{eqnarray*}
    g(x, y) = f(y) - L \cdot \mid x - y\mid \leq f(x, y) \quad \forall x \in [a, b]
\end{eqnarray*}
причём $g(y, y) = f(y)$. Из этого следует, что график функции $f(x)$ лежит выше ломаной $g(x, y)$ при всех $x \in [a, b]$ и имеет \vm{с ней общую точку} $(y, f(y))$. \vm{Данное} свойство ломаной $g(x, y)$ можно использовать для построения метода. Этот метод начинается с выбора произвольной точки $x_0 \in [a, b]$ и составления функции $g(x, x_0) = f(x_0) - L \cdot \mid x - x_0\mid = p_0(x)$. Следующая точка $x_1$ определяется из условий $p_0(x_1) = \min_{x \in [a, b]} {(p_0(x))} (x_1 \in [a, b])$, причём $x_1 = a$ или $x_1 = b$. Далее берётся новая функция $p_1(x) = \max{(g(x, x_1), p_0(x))}$, и очередная точка $x_2$ находится из условий $p_1(x_2) = \min_{x \in [a, b]} {p_1(x)} (x_2 \in [a, b])$ и т.д. 
Пусть точки $x_1, ..., x_n (n \geq 1)$ уже известны. Тогда составляется функция:
\begin{eqnarray*}
    p_n(x) = \max{(g(x, x_n), p_{n-1}(x))} = \max_{0 \leq i \leq n}{g(x, x_i)},
\end{eqnarray*}
и следующая точка $x_{n+1}$ определяется условиями:
\begin{eqnarray*}
    p_n(x_{n+1}) = \min_{x \in [a, b]}{p_n(x), \quad x_{n+1} \in [a, b]}
\end{eqnarray*}
Если минимум $p_n(x)$ достигается в нескольких точках, то в качестве $x_{n+1}$ можно взять любую из них. Т.о. метод ломаных описан. \\
Также следует упомянуть алгоритм случайного мультистарта ~\cite{AnnealingMultistart}. Случайный мультистарт - это метод глобальной оптимизации, состоящий в многократном отыскании локальных минимумов из различных начальных точек. В своем первоначальном виде он неэффективен, однако некоторые из его модификаций могут быть полезны. Основная сложность при практической реализации метода состоит в следующем: для того, чтобы с высокой надёжностью отыскать точку глобального минимума, необходимо взять количество начальных точек для локальных алгоритмов существенно больше, чем число локальных минимумов функции, которое обычно неизвестно.


\paragraph{Восстановление матрицы корреспонденци\vm{й} для г. Москвы}

В данном разделе приведены численные эксперименты для расчета матрицы корреспонденци\vm{й} для города Москвы с использованием алгоритмов, \vm{описанных} в предыдущем разделе. Данные эксперименты были проведены на основе данных, собранных в результате опроса 2013 года по Москве и Московской области. Данные представлены csv файлом с пятью полями:

\begin{enumerate}
  \item Зона i.
  \item Зона j.
  \item Число жителей i, которые ездят на работу в j.
  \item Среднее время, затраченное на поездку, в минутах.
  \item Среднее расстояние по прямой между домом и работой (в Москве точки отправления или прибытия определяются с точностью до ближайшего метро, в Зеленограде - до центра района, в области - до центра населенного пункта).
\end{enumerate}

\ai{\begin{figure}[H]
\centering
\includegraphics[width=0.7\textwidth]{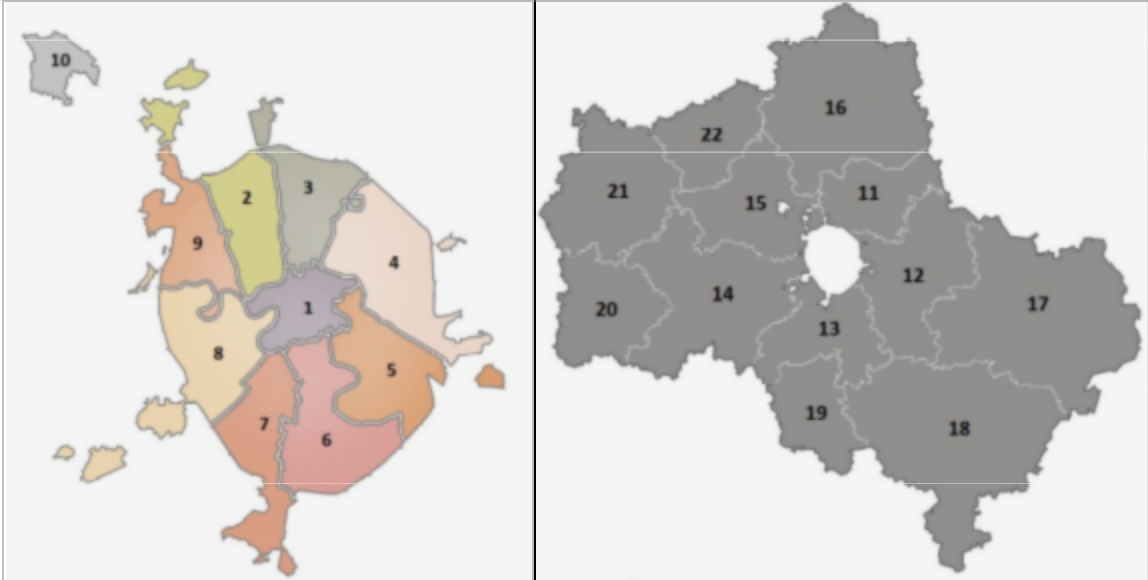}
\caption{Разбиение г.Москвы и МО по районам}
\label{fig:moscow_map}
\end{figure}
}
Отметим, что зоны выбраны достаточно большими, чтобы уменьшить случайные ошибки. В Москве (в старых границах) зоны соответствуют округам, в области нескольким средним районам. Всего есть 22 района, которые одновременно являются и источниками (\vm{место начала поездки}) и стоками (\vm{место окончания поездки}), однако не между всеми пунктами $i-j$ есть корреспонденци\vm{й}. Суммарное число участников движения $1965$.

Далее приводится описание полученных результатов для разных функций затрат.

\subparagraph{Линейная функция затрат}

Для начала в качестве функции затрат рассмотрим линейную функцию от затрат, т.е. 
$$
 T_{ij}(\alpha) = \alpha c_{ij}, 
$$
где $\alpha$- это калибруемый параметр. 

Очевидно, что чем выше альфа, тем сильнее влияют затраты на проезд по пути между источником и стоком на соответствующее значение корреспонденци\vm{й}.

Рассмотрим следующие три варианта затрат:

\begin{itemize}
    \item \textbf{Затраты - среднее время в пути.}
    
    Рассмотрим простейшую модель, где $c_{ij}$  - среднее время проезда от района $i$ до района $j$. В рамках этой модели функция затрат имеет следующий вид:
$$
 T_{ij}(\alpha) = \alpha \cdot \text{\vm{time}}_{ij}.
$$
В данной функции в качестве калибруемого параметра выступает $\alpha$. Подбор параметра происходил путем перебора по сетке $\alpha \in [0.01, 1]$ с шагом 1e-3. График зависимости невязки~\eqref{error_of_restored_matrix} от параметра $\alpha$ представлен на рис. \ref{fig:linear_cost_gamma}.

\begin{figure}[H]
    \centering
    \includegraphics[width=0.5\textwidth]{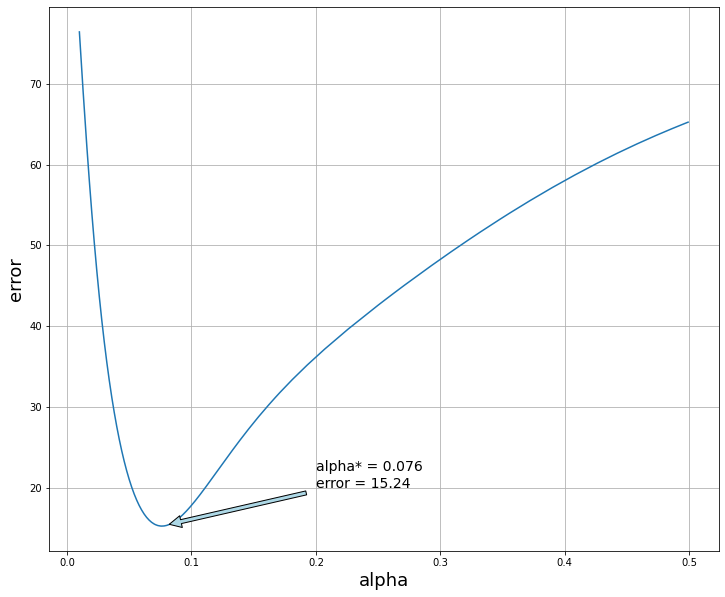}
    \caption{Значение невязки~\eqref{error_of_restored_matrix} при функции затрат $T_{ij}(\alpha) = \alpha \cdot \text{\vm{time}}_{ij}$ в зависимости от $alpha$}
    \label{fig:linear_cost_gamma}
\end{figure}

Так же в результате перебора по сетки было найдено оптимальное значение $\alpha^* = 0.076$, невязка при данном значении параметра равна $15.24$.

\item \textbf{Затраты - среднее время в пути (степенная функция затрат).}

Рассмотрим усложненную двухпараметрическую модель, где $c_{ij}$  - среднее время проезда от района $i$ до района $j$ в степении $\gamma$. В рамках этой модели функция затрат имеет следующий вид:
\begin{equation}
 T_{ij}(\alpha, \gamma) = \alpha \cdot \text{\vm{time}}^{\gamma}_{ij}.   
\label{eq:cost_alpha_gamma_time}
\end{equation}
 
В данной функции в качестве калибруемого параметра выступают два параметра: $\alpha$ и $\gamma$. Подбор параметров происходил путем перебора по сетке $\gamma \in [0.01, 1]$ с шагом 1e-2 и динамическом определении области перебора параметра $\alpha$ с целью определения окрестности, в которой достигается минимум невязки. Параметр $\gamma$ определяет ширину углубления \vm{на рисунке~\ref{fig:linear_cost_gamma}}. В таблице \ref{table:best_alpha_gamma_time_table} приведены лучшие комбинации параметров, а на графике \ref{fig:min_error_alpha_gamma_time_matrix} зависимость минимального значения невязки от $\gamma$. 

\begin{table}[H]
\caption{Сравнение невязок для разных функций затрат}
\begin{center}
\begin{tabular}{|l|l|l|}
\hline

$\gamma$ & $\alpha$ &  Невязка \\ \hline 
0.09 & 26.760 & 12.38466 \\ \hline 
0.1 & 23.770 & 12.40603 \\ \hline 
0.11 & 20.095 &12.40931 \\ \hline 
0.12 & 18.290 &12.41265 \\ \hline 
0.08 & 31.410 &12.41749 \\ \hline 
0.05 & 56.250 &12.42543 \\ \hline 
0.13 & 15.785 &12.42546 \\ \hline 
0.07 & 37.230 &12.43197 \\ \hline 
0.14 & 14.570 &12.43495 \\ \hline 
0.15 & 12.750  & 12.43975 \\ \hline 

\end{tabular}
\end{center}
\label{table:best_alpha_gamma_time_table}
\end{table}

\begin{figure}[H]
    \centering
    \includegraphics[width=0.6\textwidth]{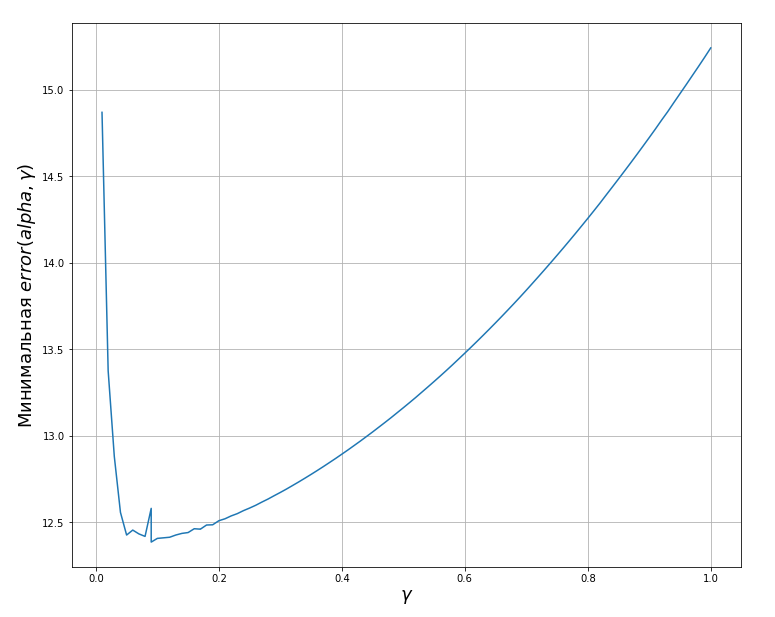}
    \caption{Оптимальное значение невязки в зависимости от $\gamma$ при функции затрат $T_{ij}(\alpha, \gamma) = \alpha \cdot \text{\vm{time}}^{\gamma}_{ij}$}
    \label{fig:min_error_alpha_gamma_time_matrix}
\end{figure}

Из графика хорошо видно, что введение степени $\gamma$ матрицы временных затрат имеет смысл для минимизации невязки. Так оптимальные значения $\gamma^* = 0.09$ и $\alpha^* = 26.760$ дают невязку равную $12.385$.

\item \textbf{Затраты - комбинация времени и расстояния.}

Введем теперь также дополнительную зависимость от расстояния в затраты. Тогда функция затрат будет иметь следующий вид:
\begin{equation}
 T_{ij}(\alpha, \beta, \gamma) = \alpha \cdot \text{\vm{time}}^{\gamma}_{ij} \cdot \text{dist}^{\beta}_{ij},  
 \label{eq:cost_alpha_gamma_beta_time_distance}
\end{equation}

где $\text{dist}_{ij}$ -- расстояние между районами $i$ и $j$, $\text{\vm{time}}_{ij}$ -- среднее время в пути между \vm{районами $i$ и $j$}.
В данной функции в качестве калибруемого параметра выступают три параметра: $\alpha$, $\beta$ и $\gamma$. Для проверки целесообразности добавления можителя $\text{dist}^{\beta}_{ij}$ будем перебирать параметр $\beta \in [0, 0.5]$ с шагом $0.001$. На графике \ref{fig:beta_distance} показано изменение невязки в зависимости от $\beta$ для оптимальных параметров при функции затрат \vm{\eqref{eq:cost_alpha_gamma_beta_time_distance}}.

\begin{figure}[H]
    \centering
    \includegraphics[width=0.6\textwidth]{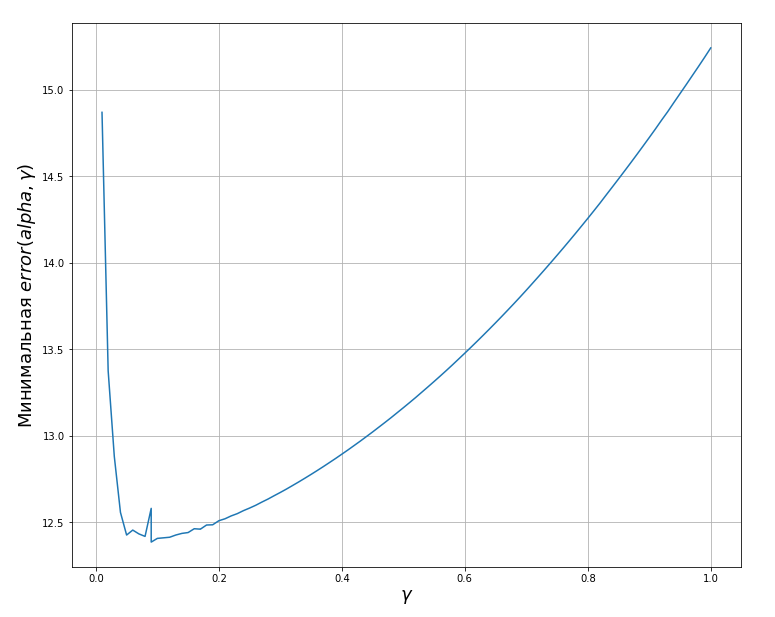}
    \caption{Изменение невязки от $\beta$ при функции затрат $T_{ij}(\alpha, \beta, \gamma) = \alpha \cdot \text{\vm{time}}^{\gamma}_{ij} \cdot \text{dist}^{\beta}_{ij}$, где $\alpha^* = 26.76$ и $\gamma^* = 0.09$}
    \label{fig:beta_distance}
\end{figure}

Видим, что добавление множителя $\text{dist}^{\beta}_{ij}$ в \eqref{eq:cost_alpha_gamma_time} имеет смысл, так как минимальная невязка достигается при $\beta = 0.005$. Сравним значения оптимальных невязок для \eqref{eq:cost_alpha_gamma_time} и \eqref{eq:cost_alpha_gamma_beta_time_distance}:

\begin{figure}[H]
    \centering
    \includegraphics[width=0.6\textwidth]{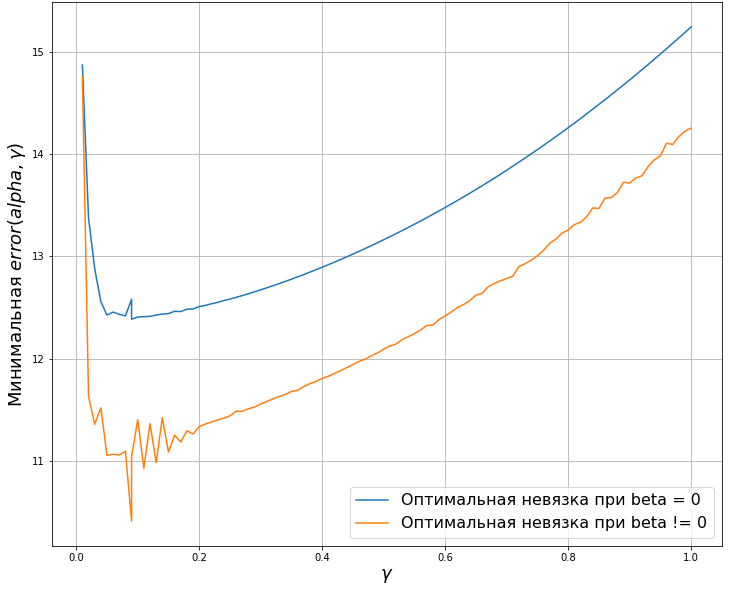}
    \caption{Сравнение оптимальная невязок  при функциях затрат $T_{ij}(\alpha, \beta, \gamma) = \alpha \cdot \text{\vm{time}}^{\gamma}_{ij} \cdot \text{dist}^{\beta}_{ij}$ и $T_{ij}(\alpha, \gamma) = \alpha \cdot \text{\vm{time}}^{\gamma}_{ij}$
    в зависимости от $\gamma$}
    \label{fig:beta_comparisons}
\end{figure}

Оптимальное значение невязки достигается при $\alpha^* = 26.76$, $\gamma* = 0.09$, $\beta^* = 0.005$ и равно 10.41226.

\end{itemize}

\subparagraph{\vm{Сумма} степенной и логарифмической функции затрат}

Следуя подходу из~\cite{gasnikov2013book} рассмотрим модель, в которой функция затрат является \vm{суммой} степенной и логарифмической, то есть 
\begin{equation}
T_{ij}(\alpha, \beta, \gamma) = \alpha c_{ij}^\gamma - \beta \ln c_{ij}. 
\label{eq:general_log}    
\end{equation}

В качестве обоснования использования данной функции можно привести следующее рассуждение: первое слагаемое отражает нежелательность больших затрат на дорогу,а второе отражает возможность найти работу на расстоянии (среднем времени в пути) порядка $c_{ij}$  от дома. Для данной функции в качестве затрат можно рассматривать затраты равные среднему времени в пути или затраты равные расстоянию между районами. Данная модель зависит от трех параметров. В случае затрат равных среднему времени в пути функция затрат будет иметь следующий вид 

\begin{equation}
T_{ij}(\alpha, \beta, \gamma) = \alpha  \text{\vm{time}}_{ij}^\gamma - \beta \ln  \text{\vm{time}}_{ij}. 
\label{eq:cost_alpha_gamma_beta_time}
\end{equation}

\begin{figure}[H]
    \centering
    \includegraphics[width=0.6\textwidth]{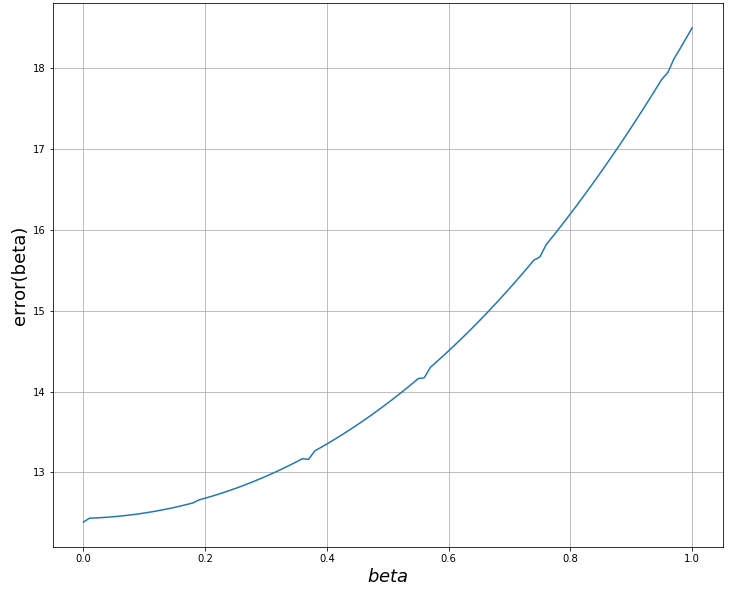}
    \caption{Невязка в зависимости от $\beta$ при функции затрат $T_{ij}(\alpha, \gamma, \beta) = \alpha \cdot \text{\vm{time}}^{\gamma}_{ij} - \beta \ln  \text{\vm{time}}_{ij}$ при оптимальных $(\alpha, \gamma)$}
    \label{fig:error_beta_log_time}
\end{figure}

По результатам проведенных экспериментов можно заключить, что оптимальные значения невязки достигаются при $\beta = 0$, то есть, когда функция затрат \eqref{eq:cost_alpha_gamma_beta_time} вырождается в \eqref{eq:cost_alpha_gamma_time}. На графике \ref{fig:error_beta_log_time} показано, как ухудшается невязка при росте $\beta$ при оптимальных $(\alpha, \gamma)$.

Если в \eqref{eq:general_log} $c_{ij} = \text{\vm{dist}}_{ij}$, то получим функцию затрат от среднего расстояния в пути:

\begin{equation}
T_{ij}(\alpha, \beta, \gamma) = \alpha  \text{\vm{dist}}_{ij}^\gamma - \beta \ln  \text{\vm{dist}}_{ij}. 
\label{eq:cost_alpha_gamma_beta_dist}
\end{equation}
Перебор по параметрам $\alpha \in [0.01, 10], \beta \in [0, 0.5], \gamma \in [0, 1]$ выдал оптимальную точку $(\alpha^* = 3.01, \beta^* = 0, \gamma^* = 0.25)$ с невязкой 4.84729.



\subparagraph{Сравнение функций затрат}

В данном разделе приведем сравнение результатов, полученных выше. Обозначим $t_{ij}$- среднее время в пути между районами $i$ и $j$, а $d_{ij}$- среднее расстояние между районами $i$ и $j$. В таблице~\ref{error_table} приведены минимальные значения невязок для различных функций  затрат и соответствующие этим значениям невязок оптимальные параметры.  

\begin{table}[H]
\caption{Сравнение невязок для разных функций затрат.}
\begin{center}
\label{error_table}
\begin{tabular}{|l|l|l|}
\hline
Функция затрат $T_{ij}(\alpha, \beta, \gamma)$      & Невязка & Оптимальные параметры                  \\ \hline
$\alpha \cdot t_{ij}$                               &$15.24$    & $\alpha =0.076$                             \\ \hline
$\alpha \cdot t_{ij}^{\gamma}$                       & $12.38$   & $\alpha = 26.76, \, \gamma = 0.09$               \\ \hline
$\alpha \cdot t_{ij}^{\gamma} \cdot d_{ij}^{\beta}$ & 10.41   & $\alpha = 26.76, \, \beta = 0.005, \,  \gamma = 0.09$ \\ \hline
$\alpha t_{ij}^\gamma - \beta \ln t_{ij}$           & 12.38   & $\alpha = 26.76, \, \beta = 0, \,  \gamma = 0.09$ \\ \hline
$\alpha d_{ij}^\gamma - \beta \ln d_{ij}$           & 4.84 & $\alpha = 3.01, \, \beta = 0, \,  \gamma = 0.25$ \\ \hline
\end{tabular}
\end{center}
\end{table}

По итогам проведенных экспериментов заключаем, что оптимальной функцией затрат для восстановления матрицы корреспонденций по имеющимся данных, является $T_{ij} = 3.01 \cdot d_{ij}^{0.25}$. Заметим, что функция затрат не зависит от времени. Это может говорить как о недостаточности объема временных данных, все-таки для измерения средних значений по времени требуется больше наблюдений, так и о необходимости продолжить исследования по подбору иных функций затрат.    

\paragraph{Заключение}
В статье была предложена и научно обоснована технология калибровки моделей расчета матрицы корреспонденций. В отличие от известных ранее (однопараметрических) моделей расчета матрицы корреспонден\vm{ций} в статье рассматривались многопараметрические модели. Также отметим важное отличие предлагаемых исследований от известных ранее. Проводимые ранее исследования в данном направлении, как правило, не обосновывались строго. Более того, в подавляющем числе случаев это считалось простой задачей и про это не особо подробно писали в соответствующих работах. Однако, уже для случая двух параметров задача становится существенно более затратной и от выбора способа ее численного решения существенно зависит эффективность подхода. Собственно, вычислительным аспектам в данной статье и было уделено наибольшее внимание. Отметим, что проводимые в данной работе исследования (точнее их теоретическое обоснование) стало возможным благодаря современным достижениям, как в области математического моделирования транспортных потоков (появился новый цикл многопараметрических моделей расчета матрицы корреспонден\vm{ций}, обобщающих энтропийную модель), так и тонкому анализу скорости сходимости алгоритм Синхорна, использующегося для расчета матрицы корреспонденций (напомним, что модель расчета матрицы корреспонденций, сводится к системе 2n нелинейных уравнений, для решения которой и используется алгоритм Синхорна).

\dotfill

Задача была предложена Е.А. Нурминским и А.В. Гасниковым.

\end{document}